\documentclass[10pt]{article}

\usepackage{amsmath,amssymb,amsthm}
\usepackage{amsxtra}
\usepackage{amsfonts}
\usepackage{amssymb}
\usepackage{url}
\usepackage{hyperref}
\usepackage{esint}
\usepackage{eucal}
\usepackage{mathrsfs}
\usepackage{cite}
\usepackage{graphicx,graphics}
\newcommand\fig[1] {{\rm Figure}~\ref{fig:#1}}
\newcommand\labfig[1] {\label{fig:#1}}

\newcommand\eq[1] {(\ref{#1})}
\newcommand{\bfm}[1]{\mbox{\boldmath ${#1}$}}

\newtheorem{theorem}{Theorem}[section]

\newcommand{\nonum}{\nonumber \\}
\newcommand{\beqa}{\begin{eqnarray}}
\newcommand{\eeqa}[1]{\label{#1}\end{eqnarray}}
\newcommand{\beq}{\begin{equation}}
\newcommand{\eeq}[1]{\label{#1}\end{equation}}
\newcommand{\Grad}{\nabla}
\newcommand{\Div}{\nabla \cdot}

\newcommand{\Tr}{\mathop{\rm Tr}\nolimits}

\newcommand{\Md}{\partial}

\newcommand{\Ga}{\alpha}
\newcommand{\Gb}{\beta}
\newcommand{\Gd}{\delta}

\newcommand{\Gg}{\gamma}
\newcommand{\Gc}{\chi}

\newcommand{\Gk}{\kappa}

\newcommand{\Gl}{\lambda}
\newcommand{\Gn}{\eta}
\newcommand{\Gm}{\mu}

\newcommand{\Gs}{\sigma}

\newcommand{\GT}{\Theta}

\newcommand{\GO}{\Omega}

\newcommand{\BGe}{\bfm\epsilon}

\newcommand{\BGn}{\bfm\eta}

\newcommand{\BGs}{\bfm\sigma}

\newcommand{\CE}{{\cal E}}

\newcommand{\CQ}{{\cal Q}}

\newcommand{\CS}{{\cal S}}

\newcommand{\bpm}{\begin{pmatrix}}
\newcommand{\epm}{\end{pmatrix}}

\def\Ba{{\bf a}}
\def\Bb{{\bf b}}

\def\Be{{\bf e}}

\def\Bh{{\bf h}}

\def\Bk{{\bf k}}

\def\Bn{{\bf n}}

\def\Bq{{\bf q}}

\def\Bt{{\bf t}}
\def\Bu{{\bf u}}

\def\Bx{{\bf x}}

\def\BA{{\bf A}}
\def\BB{{\bf B}}
\def\BC{{\bf C}}

\def\BI{{\bf I}}

\def\BK{{\bf K}}

\def\BS{{\bf S}}

\def\half{{\scriptstyle{1\over 2}}}

\newcommand{\ajj}[1]{\textcolor{black}{#1}}
\usepackage{graphics}
\usepackage{amssymb}
\usepackage{stackengine}
\usepackage{color}
\setstackgap{S}{1pt}
\title{Near optimal pentamodes as a tool for guiding stress while minimizing compliance in $3d$-printed materials: a complete solution to the weak $G$-closure problem for $3d$-printed materials.}
\date{}
\begin{document}
\maketitle
\vskip -.5cm
\centerline{\large
Graeme W. Milton \footnote{Department of Mathematics, University of Utah, USA -- milton@math.utah.edu,}
\,\, and \,\,
Mohamed Camar-Eddine\footnote{Institut de Recherche Math\'ematique de Rennes, INSA de Rennes, FRANCE -- Mohamed.Camar-Eddine@insa-rennes.fr.}}
\vskip 1.cm
\begin{abstract}
For a composite containing one isotropic elastic material, with positive Lame moduli, and void, with the elastic material occupying a prescribed volume fraction $f$,
and with the composite being subject to an average stress, $\BGs^0$, Gibiansky, Cherkaev, and Allaire provided a sharp lower bound $W_f(\BGs^0)$ on the 
minimum compliance energy $\BGs^0:\BGe^0$, in which $\BGe^0$ is the average strain. Here we show these bounds also provide sharp bounds on the
possible $(\BGs^0,\BGe^0)$-pairs that can coexist in such composites, and thus solve the weak $G$-closure problem for $3d$-printed materials. The
materials we use to achieve the extremal $(\BGs^0,\BGe^0)$-pairs are denoted as near optimal pentamodes.
We also consider two-phase composites containing this isotropic elasticity material and a rigid phase with the elastic material occupying a prescribed volume fraction $f$,
and with the composite being subject to an average strain, $\BGe^0$. For such composites, Allaire and Kohn provided a sharp lower bound $\widetilde{W}_f(\BGe^0)$ on the 
minimum elastic energy $\BGs^0:\BGe^0$. We show that these bounds also provide sharp bounds on the
possible $(\BGs^0,\BGe^0)$-pairs that can coexist in such composites of the elastic and rigid phases, and thus solve the weak $G$-closure problem in this case too.
The materials we use to achieve these extremal $(\BGs^0,\BGe^0)$-pairs are denoted as near optimal unimodes.
\end{abstract}
%%%%%%%%%%%%%%%%%%%%%%%%%%%%%%%%%%%%%%%%%%%%%%%%%%%%%%%%%%%%%%%%%%%%%%%%
\section{Introduction}
\setcounter{equation}{0}
%%%%%%%%%%%%%%%%%%%%%%%%%%%%%%%%%%%%%%%%%%%%%%%%%%%%%%%%%%%%%%%%%%%%%%%%%%%%%%%%%%%%%%%%%%%%%%%%%%%%%%%%%%%
Consider the set of all $3d$-printed periodic materials built from a single elastically isotropic material 
with elasticity tensor $\BC_1$ occupying a volume fraction $f$ with void occupying the remaining volume fraction $1-f$, which is known as the porosity.
Each material in this set has an effective elasticity tensor $\BC^*$ and a natural question
to ask is: what are the possible (average stress, average strain)-pairs that can exist in the composite allowing for all possible microgeometries
having the volume fraction $f$? By average we mean an average over the unit cell of periodicity, and it is assumed that the macroscopic
variation of the local stress $\BGe(\Bx)$ and strain $\BGs(\Bx)$ is such that these average quantities vary very slowly from cell to neighboring cells.
To compute the possible (average stress, average strain)-pairs it suffices to assume that  $\BGe(\Bx)$ and $\BGs(\Bx)$ are periodic, and
then their averages $\BGe^0$ and $\BGs^0$ are the same for any cell. By definition of the effective elasticity tensor $\BC^*$, one has
$\BGs^0=\BC^*\BGe^0$. The problem of characterizing the set of possible average field pairs has been called the weak $G$-closure problem by Cherkaev \cite{Cherkaev:2000:VMS}. 
Our assumption that the geometry is periodic is made for ease of analysis, and does not restrict the set of possible effective elasticity tensors and hence does not restrict the 
set of possible (average stress, average strain)-pairs \cite{Raitums:2001:LRC,Allaire:2002:SOH}.

The upshot of this paper is that it completely solves the weak $G$-closure elasticity problem for $3d$-printed materials.
For two-dimensional porous media (or three-dimensional porous media having microstructure independent of one coordinate) the weak $G$-closure problem
was solved in \cite{Milton:2003:RAS}, giving an essentially complete characterization of
the set of possible (average stress, average strain)-pairs. We mention too that bounds have been derived
in the case where the second phase is an elastic material, rather than void \cite{Milton:2011:BVF}, and progress has also been made
on the problem of bounding average field pairs for nonlinear composites 
\cite{Milton:2000:BCN,Talbot:2004:BEC,Peigney:2005:PBM,Peigney:2008:RSC,Bhattacharya:2010:SIP}.

The question can be reinterpreted as: Given a prescribed average stress $\BGs^0$, what is the range $\CE_f(\BGs^0)$ of values of $\BGe^0$ such that
there exists a $3d$-printed periodic material having the desired volume fraction $f$ with an effective elasticity tensor $\BC^*$
such that $\BGs^0=\BC^*\BGe^0$? \ajj{This weak $G$-closure elasticity problem is considerably simpler than the $G$-closure problem of trying to find the
set of possible tensors $\BC^*$. Anisotropic elasticity tensors in $3d$ are characterized by
18 invariants. Thus assuming the volume fraction and the elastic moduli of the phases are given, the $G$-closure is represented by a set in an 18-dimensional
space. For the weak $G$-closure one can assume without loss of generality that $\BGs^0$ is diagonal and normalized so $\Tr{\BGs^0\BGs^0}=1$. This leaves two
parameters free. Then the $3\times 3$ matrix $\BGe^0$ has 6 elements. Thus, the set of (average stress, average strain)-pairs can be represented as a set
in an 8-dimensional space, a considerable simplification over the 18-dimensional space needed for the $G$-closure.}

One constraint on the set $\CE_f(\BGs^0)$ is immediately apparent: the known elastic compliance energy bounds of
Gibiansky and Cherkaev \cite{Gibiansky:1987:MCE}, and Allaire \cite{Allaire:1994:ELP} should be satisfied:
\beq W_f(\BGs^0)\leq\BGs^0:\BGe^0, \eeq{0.1}
where their explicit formula for the function $W_f(\BGs^0)$ will be given in the next section. Here, and afterwards, the colon $:$ represents a double
contraction of indices. Thus $\BA:\BB$ represents the natural inner product of the space of symmetric matrices with $\BA:\BB=\Tr(\BA\BB)$.
Allaire's result was a specialization of a more general result of Allaire and Kohn \cite{Allaire:1993:OBE}, applicable
to two-phase composites where both phases are elastic: see their equation (6.9). 

The bound \eq{0.1} is a linear constraint on $\BGe^0$ that
constrains $\BGe^0$ to lie on one side of a hyperplane in the six-dimensional space of symmetric $3\times 3$ matrices.
It is the objective of this paper to show
that the bound is asymptotically sharp, and thus essentially completely characterizes $\CE_f(\BGs^0)$, in the sense that given any
 $\BGs^0$ and any $\BGe^0$ such that \eq{0.1} holds, there exists a sequence of microstructures
having effective tensors $\BC^*_\Gd$ and porosities $1-f$ such that $\BC^*_\Gd\BGe^0$ converges to $\BGs^0$,
as $\Gd\to 0$. The materials we will use fall in the class of pentamode materials, and if they are such that \eq{0.1} is almost
satisfied as an equality, they will be called near optimal pentamode materials. A blueprint for the theoretical design of near
optimal pentamode materials is implicit in the constructions given in \cite{Milton:2016:PEE}. Basically, one takes the homogenized
sequential laminate attaining the bounds \eq{0.1} and inserts into it sets of parallel walls with thickness much smaller than their separation.
The walls themselves contain a suitably chosen pentamode material that supports the desired stress $\BGs^0$ but otherwise allows slip.
In the end one obtains a sort of cellular material where the interior of the cells contains the homogenized sequential laminate attaining the energy
bounds \eq{0.1} and the cell walls are thin and contain a suitably constructed pentamode material. 

\ajj{In the specific case where
$\BGs^0$ is proportional to the identity tensor the optimal pentamodes are elastically isotropic materials that (asymptotically, as $\Gd\to 0$)
attain the upper Hashin-Shtrikman bulk modulus bounds \cite{Hashin:1963:VAT} yet have shear modulus approaching zero. Following the prescription laid out in \cite{Milton:2016:PEE},
they are obtained by taking a homogenized elastically isotropic 
material that attains the upper Hashin-Shtrikman bulk modulus bound, and inserting in it systems of thin parallel walls with a sufficient number of
orientations to ensure the elastic isotropy of the resulting composite. The walls themselves contain a microstructure of parallel cylinders of phase 1
embedded in phase 2, with the cylinder axes being parallel to the normal $\Bn$ to the walls, so that the walls support uniaxial compression but not shear
(such walls, not necessarily thin, also appear in the constructions of Sigmund \cite{Sigmund:2000:NCE}).
Of course, by mixing this optimal pentamode with a homogenized material simultaneously attaining both the bulk and shear modulus upper bounds one
obtains elastically isotropic materials with bulk modulus (asymptotically, as $\Gd\to 0$) attaining the upper Hashin-Shtrikman bulk modulus bound, and
with any desired shear modulus between zero and the upper Hashin-Shtrikman shear modulus bound. The possible bulk and shear moduli of elastically isotropic
$2d$ and $3d$-printed elastic materials have also been researched in \cite{Milton:1992:CMP,Huang:2011:TDM,Andreassen:2014:DME,Xia:2015:DMT,Berger:2017:MMT,Ostanin:2017:PCC} (see
also references therein, and for results in the case where the second phase is not void see also, for example, 
\cite{Hashin:1963:VAT,Berryman:1988:MRC,Milton:1992:CMP,Cherkaev:1993:CEB,Sigmund:1995:TMP,Sigmund:2000:NCE}).}

The class of pentamode materials, materials that can support any desired average stress loading $\BGs^0$, but which are easily compliant
to any other loadings, were first introduced in \cite{Milton:1995:WET} as a means to constructing materials with any desired positive definite elasticity
tensor, built from one sufficiently compliant phase and one sufficiently stiff phase. Their effective elasticity tensor $\BC^*$ has 5 small
eigenvalues (hence the name pentamode) and they include materials that can only support
hydrostatic loads, an example of which was obtained independently by Sigmund \cite{Sigmund:1995:TMP}. Liquids and gels are also examples of
materials that can essentially only support hydrostatic loadings. \ajj{Another elementary example of a pentamode material is a periodic 
array of cylinders of the elastic material, with axis parallel to a vector $\Bt$: this can support stresses proportional to the rank one tensor $\Bt\otimes\Bt$
(here, and later, $\otimes$ denotes the tensor product).}
So the most interesting pentamodes are those that support anisotropic loadings \ajj{that are not rank 1}. 
It was recognized that pentamodes might  be useful in cloaking  \cite{Norris:2008:ACT}, 
and pentamodes were first built using a three dimensional lithography technique \cite{Kadic:2012:PPM} leading to an ``unfeelability'' cloak \cite{Buckmann:2014:EMU}.
Other investigations of pentamode materials and their applications include, e.g., \cite{Scandrett:2010:ACU,Scandrett:2011:BOP,Cipolla:2011:DIP,Martin:2012:PBS,Gokhale:2012:STP,Kadic:2013:AVT,Schittny:2013:EMM,Layman:2013:HAE,Mejica:2013:CSP,Buckmann:2014:EMU,Kadic:2014:PPM,Cai:2015:PMA,Chen:2015:LPA,Tian:2015:BMA,Fabbrocino:2016:SAP,Amendola:2016:ERA,Huang:2016:PPA,Hedayati:2017:AMM}.

We also consider mixtures of the elastically isotropic phase 1 with elasticity tensor $\BC_1$ with an almost rigid phase 2 with elasticity tensor $\BC_2=\Gd\BC_0$, and
ask: in the limit $\Gd\to\infty$ what are the set of possible $(\BGs^0, \BGe^0)$ pairs when the material with tensor $\BC_1$ occupies the volume fraction $f$? 
Specifically, considering $f$ and $\BGe^0$ to be fixed, we find the range $\CS_f(\BGe^0)$ of possible values of $\BGs^0=\BC^*\BGe^0$
as the microgeometry ranges over all configurations in which phase 1 occupies the volume fraction $f$.
The known elastic energy bounds of Allaire and Kohn \cite{Allaire:1993:OBE} provide an inequality of the form
\beq \widetilde{W}_f(\BGe^0)\leq\BGs^0:\BGe^0, \eeq{0.10}
where their formula for $\widetilde{W}_f(\BGe^0)$ will be given in Section 3. We will 
show this bound is asymptotically sharp, and thus essentially completely characterizes $\CS_f(\BGs^0)$, in the sense that given any
$\BGs^0$ and any $\BGe^0$ such that \eq{0.1} holds, there exists a sequence of microstructures
having effective tensors $\BC^*_\Gd$ and with phase 1 occupying a volume fraction $f_\Gd$ such that $\BC^*_\Gd\BGe^0$ converges to $\BGs^0$, and $f_\Gd$ converges to $f$
as $\Gd\to \infty$. The materials we will use fall in the class of unimode materials, and if they are such that \eq{0.10} is almost
satisfied as an equality, they will be called near optimal unimode materials. A blueprint for the theoretical design of near
optimal unimode materials is implicit in the constructions given in \cite{Milton:2016:TCC}. Basically, one takes the homogenized
sequential laminate attaining the bounds \eq{0.1} and inserts into it sets of parallel walls with thickness much smaller than their separation.
The walls themselves contain a suitably chosen unimode material that can slip, but only with $\BGe^0$ as the associated strain in the walls.
In the end one obtains a sort of cellular material where the interior of the cells contains the homogenized sequential laminate attaining the energy
bounds \eq{0.10} and the cell walls are thin and contain a suitably constructed unimode material.

Unimode materials have the property that there is only one easy mode of deformation: the effective compliance tensor $\BS^*=(\BC^*)^{-1}$ has only
one large eigenvalue. Examples of them (without sliding surfaces) were first given in \cite{Milton:1992:CMP}. Later it was shown that there are
unimodes that allow any desired strain $\BGe^0$ as their easy mode of deformation \cite{Milton:1995:WET}. One can generalize the idea
of unimode materials to finite deformations, and in that context they are materials whose deformation  is confined to single 
one-dimensional trajectory in the space of Cauchy-Green tensors. In fact one can build up an entire algebra of their macroscopic responses, including addition, subtraction, and composition, 
and thereby show that any conceivable deformation trajectory (in the set of positive definite Cauchy-Green tensors) is realizable to an arbitrarily high degree of approximation
\cite{Milton:2013:CCMa}.

\ajj{In the specific case where
$\BGe^0$ is proportional to the identity tensor the optimal unimodes are elastically isotropic materials that (asymptotically, as $\Gd\to \infty$)
attain the lower Hashin-Shtrikman bulk modulus bounds \cite{Hashin:1963:VAT} yet have shear modulus approaching infinity. Following the prescription laid out in \cite{Milton:2016:PEE},
they are obtained by taking a homogenized elastically isotropic 
material that attains the lower Hashin-Shtrikman bulk modulus bound, and inserting in it systems of thin parallel walls with a sufficient number of
orientations to ensure the elastic isotropy of the resulting composite. The walls themselves contain a microstructure of unimode material
so that the walls allow uniaxial compression but are resistant to shear.
Of course, by mixing this optimal unimode with a homogenized material simultaneously attaining both the bulk and shear modulus lower bounds one
obtains elastically isotropic materials with bulk modulus (asymptotically, as $\Gd\to 0$) attaining the lower Hashin-Shtrikman bulk modulus bound, and
with any desired shear modulus between infinity and the lower Hashin-Shtrikman shear modulus bound.}

%%%%%%%%%%%%%%%%%%%%%%%%%%%%%%%%%%%%%%%%%%%%%%%%%%%%%%%%%%%%%%%%%%%%%%%%
\section{The relevance to a wide class of optimization questions}
\setcounter{equation}{0}
%%%%%%%%%%%%%%%%%%%%%%%%%%%%%%%%%%%%%%%%%%%%%%%%%%%%%%%%%%%%%%%%%%%%%%%%%%%%%%%%%%%%%%%%%%%%%%%%%%%%%%%%%%%

This section can be skipped by those readers not interested in structural optimization, but it does give insight into the 
importance of trying to characterize the set $\CE_f(\BGs^0)$ for all $3\times 3$ symmetric matrices $\BGs^0$ and all
volume fractions $f$ between 0 and 1. \ajj{Rather than solving specific elasticity examples the purpose of the section is just to indicate
the breadth of optimization problems that may be solved once one knows the elastic weak $G$-closure.} 
As motivation, we begin by looking at the simpler thermal conductivity problem (or the equivalent electrical conductivity problem) 
for composites of two isotropic materials with thermal conductivities $k_1$ and $k_2$, mixed in proportions $f$ and $1-f$.
The equations take the form
\beq \Bq(\Bx)=k(\Bx)\Be(\Bx),\quad \Div \Bq=0, \quad \Be=-\Grad T, \eeq{0.1a}
where $T$ is the temperature field, with gradient $\Be(\Bx)$, $\Bq(\Bx)$ is the heat current, and $k(\Bx)$ is $k_1$ in phase 1 and $k_2$ in phase 2.
Then an analogous question is the following: supposing $\Be(\Bx)$ and $\Bq(\Bx)$ are periodic and that
the average $\Be^0$ of $\Be(\Bx)$ is prescribed, what is the range of values $\CQ_f(\Be^0)$ that the average $\Bq^0$ of $\Bq(\Bx)$
takes as the geometry varies over all periodic microstructures with phase 1 having the required volume fraction $f$? The answer, as illustrated
in \fig{1},
was provided by Ra{\u\i}tum \cite{Raitum:1978:EEP,Raitum:1983:QES} who showed that $\Bq^0\in\CQ_f(\Be^0)$ if and only if
\beq |\Bq^0-\half(\Gk^++ k^-)\Be^0|^2\leq ( k^+- k^-)^2|\Be^0|^2/4,
\eeq{0.2}
where 
\beq  k^+=f k_1+(1-f) k_2,\quad  k^-=(f/ k_1+(1-f)/ k_2)^{-1},
\eeq{0.3}
and for any $\Bq^0$ and $\Be^0$ such that equality holds in \eq{0.2} one can find a simple laminate of the two phases with effective thermal conductivity
tensor $\BK^*$ having one eigenvalue taking the value $k^-$ and remaining (one or two according to whether the dimension $d$ is two or three) eigenvalues taking
the value $k^+$, such that $\Bq^0=\BK^*\Be^0$. See also the related results in \cite{Murat:1985:CVH,Gibiansky:1988:OFH,Tartar:1987:AOO,Tartar:1994:ROD,Tartar:1995:RHM},
sections 3.24, 5.2-5.4 of \cite{Cherkaev:2000:VMS} and section 22.4 of \cite{Milton:2002:TOC}. \ajj{Note that the classic Wiener bound \cite{Wiener:1912:TMF} implies 
$k^-|\Be_0|^2\leq\Bq_0\cdot\Be_0\leq k^+|\Be_0|^2$ and confines $\Bq_0$ to lie between two planes, not the disk or sphere implied by \eq{0.2} (see \fig{1}).
In fact, though, additional arguments do show that \eq{0.2} follows from the Wiener bounds (see, for example, section 22.4 of \cite{Milton:2002:TOC})}. 

\begin{figure}[!ht]
\centering
\includegraphics[width=0.5\textwidth]{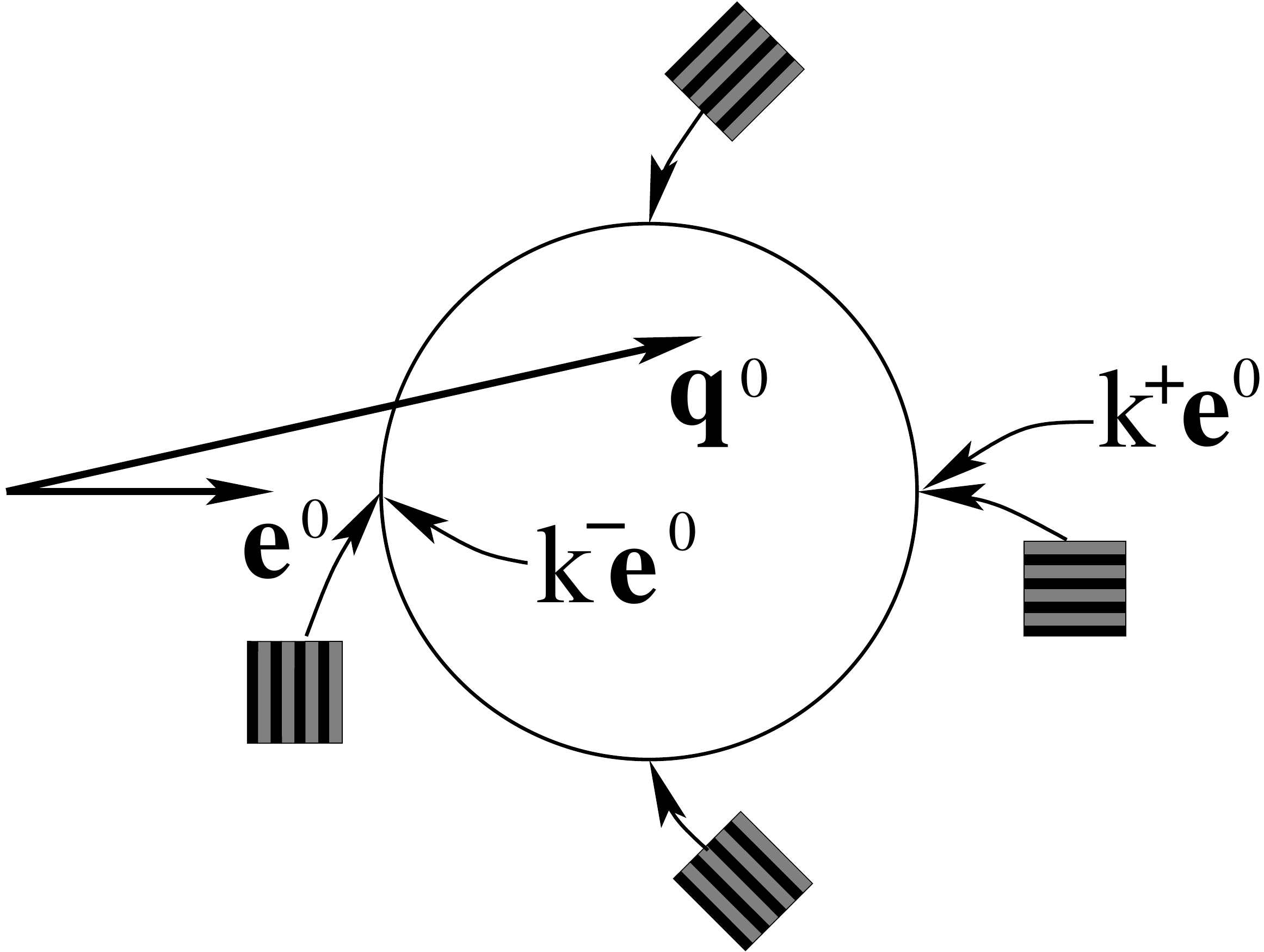}
\caption{The bounds on (average heat current, average temperature gradient) pairs $(\Bq^0,\Be^0)$ are shown.  In two-dimensions, for a given $\Be^0$, the average current $\Bq^0$ 
takes values ranging throughout the disk (or sphere in three-dimensions) as the periodic microgeometry is varied while
keeping the volume fraction $f$ of phase 1 fixed. Values on the boundary are achieved
when the geometry is a simple laminate, suitably oriented.} 
\labfig{1}
\end{figure}

\ajj{As a very simple, but illuminating, example suppose that phase 2 is a thermally insulating material, with $k_2=0$. Then the bounds \eq{0.2} imply that for a fixed heat flux $\Bq^0$ the temperature gradient $\Be^0$ is
constrained by
\beq \Bq^0\cdot\Be^0\geq (fk_1)^{-1}|\Bq^0|^2.
\eeq{0.4}
Similar to our problem, this is a linear constraint on $\Be^0$. Now suppose that we desire to guide the macroscopic heat current $\Bq(\Bx)$ parallel to a vector field $\Ba(\Bx)$.
The simplest way to do this is to choose the microgeometry so that the local effective thermal conductivity $\Bk^*(\Bx)$ takes the form
\beq \Bk^*(\Bx)=\Ga(\Bx)\Ba(\Bx)\otimes\Ba(\Bx). \eeq{0.3a}
Then $\Bq(\Bx)=-\Ga(\Bx)[\Ba(\Bx)\cdot\Grad T(\Bx)]\Ba(\Bx)$ is indeed parallel to $\Ba(\Bx)$ and we are left solving the equation
\beq 0=-\Div\Bq=\Div\{\Ga(\Bx)[\Ba(\Bx)\cdot\Grad T(\Bx)]\Ba(\Bx)\}, \eeq{0.3b}
for $T(\Bx)$. The desired effective thermal conductivity $\Bk^*$ given by \eq{0.3a} is achieved by a laminate structure in $2d$ or in $3d$ by
a structure consisting of parallel wires of phase 1 in a matrix of the thermally insulating phase 2. Of course, generally we want
to guide the heat current with minimal thermal resistance, and this can be done with a microgeometry attaining the Wiener bound
$\Bq(\Bx)\cdot\Be(\Bx)=(fk_1)^{-1}|\Bq(\Bx)|^2$ where we keep track of the local volume fraction $f(\Bx)$ as we may additionally like to fix the
overall proportions of the two components as the overall weight, or cost, may be important assuming the two components have different
weights or costs.}

\ajj{The situation with guiding stress is very similar. One can think of the three columns of the macroscopic stress $\BGs(\Bx)$ as three current fields that we would
like to guide, while keeping the symmetry of $\BGs(\Bx)$. One could specify a desired symmetric valued matrix field $\BA(\Bx)$ and demand
that $\BGs(\Bx)=\Gb(\Bx)\BA(\Bx)$ everywhere, with $\Gb(\Bx)$ being the local proportionality factor. As pointed out by Norris \cite{Norris:2008:ACT}
this may not be possible as the three components of $\Div\BGs$ that must be zero overconstrain
the single free scalar parameter $\Gb(\Bx)$. So instead one can require that $\BA(\Bx)$ satisfy $\Div\BA=0$, and that $\Gb(\Bx)$ is a constant $\Gb$.
Analogously to \eq{0.3a} we can get the desired macroscopic stress field if the local microgeometry is a pentamode material with elasticity tensor
\beq \BC^*(\Bx)=\Ga(\Bx)\BA(\Bx)\otimes\BA(\Bx). \eeq{0.3c}
Then the macroscopic displacement field (which, like in a fluid, is not uniquely determined) needs to be such that
\beq \Ga(\Bx)\Tr[\BA(\Bx)\Grad\Bu(\Bx)]=\Gb. \eeq{0.3d}
The quantity $\Gb$ is then the analog of pressure in a fluid.
Now the analogous goal to minimizing thermal resistance is to maximize the overall stiffness, in the sense of minimizing 
the local compliance energy $\BGs(\Bx):\BGe(\Bx)$ for a given macroscopic stress field $\BGs(\Bx)$ and given volume fraction. 
That is why it is advantageous to use optimal pentamodes in the
construction as they achieve the lowest possible values of $\BGs(\Bx):\BGe(\Bx)$, namely $\widetilde{W}_f(\BGs(\Bx))$.}

As a further example of the use of the bounds \eq{0.2}, consider the following two-dimensional thermal shielding problem shown in \fig{2}. 
After a suitable choice of spatial coordinates, one considers a rectangle whose left and right hand sides at $x_1=0$ and $x_1=2w$ have a uniform flux of heat across it
from left to right: $q_1(\Bx)=1$ for all $\Bx=(0,x_2)$ and $\Bx=(2w,x_2)$, $x_2\in(-1,1)$. The top and bottom sides of the rectangle at $x_2=1$ and $x_2=-1$, respectively, are thermally
insulating, so that there is no flux of current across those boundaries: $q_2(\Bx)=0$ for $\Bx=(x_1,1)$ and $\Bx=(x_1,-1)$ with $x_1\in (0,2w)$.
Inside the rectangle the thermal conductivity equations \eq{0.1a}
hold, where $k(x)$ takes the value $k_1$ in phase 1 and $k_2$ in phase 2. The objective is to configure the phases so as to
shield the interval $x_1=w$ and $x_2\in[-a,a]$ from the heat current, i.e, minimize the total heat flux across this
interval. Assuming reflection symmetry of the solution about the lines $x_1=w$ and $x_2=0$, this is the same as maximizing the total heat flux 
across the interval $x_1=w$ and $x_2\in[a,1]$ with no flux of current across the line $x_2=0$ and with the line $x_1=w$ being an isotherm,
at temperature $T=0$ with a suitable choice of temperature scale. That problem has been solved by 
Gibiansky, Lurie, and Cherkaev \cite{Gibiansky:1988:OFH} (see also Sect.6.4 in \cite{Lurie:1997:ECC} and Sect.5.4 in \cite{Cherkaev:2000:VMS}).
A rough approximation to the ideal geometry is shown in \fig{2}.
The bottom line is that composites comprised of laminates of two phases are needed at the microscale in the construction
precisely because they are the geometries attaining the bounds \eq{0.2}. Laminate geometries are the best microstructures for
guiding current. We remark that while this solution provides an optimal thermal shield, it is only optimized for one set of boundary 
conditions. It is also different to the question of thermal cloaking, where one has the additional objective of minimizing
the perturbation to the downstream heat flux \cite{Greenleaf:2003:ACC,Greenleaf:2003:NCI,Liu:2010:NSA,Guenneau:2012:TTC,Schittny:2013:ETT}.

\begin{figure}[!ht]
\centering
\includegraphics[width=0.5\textwidth]{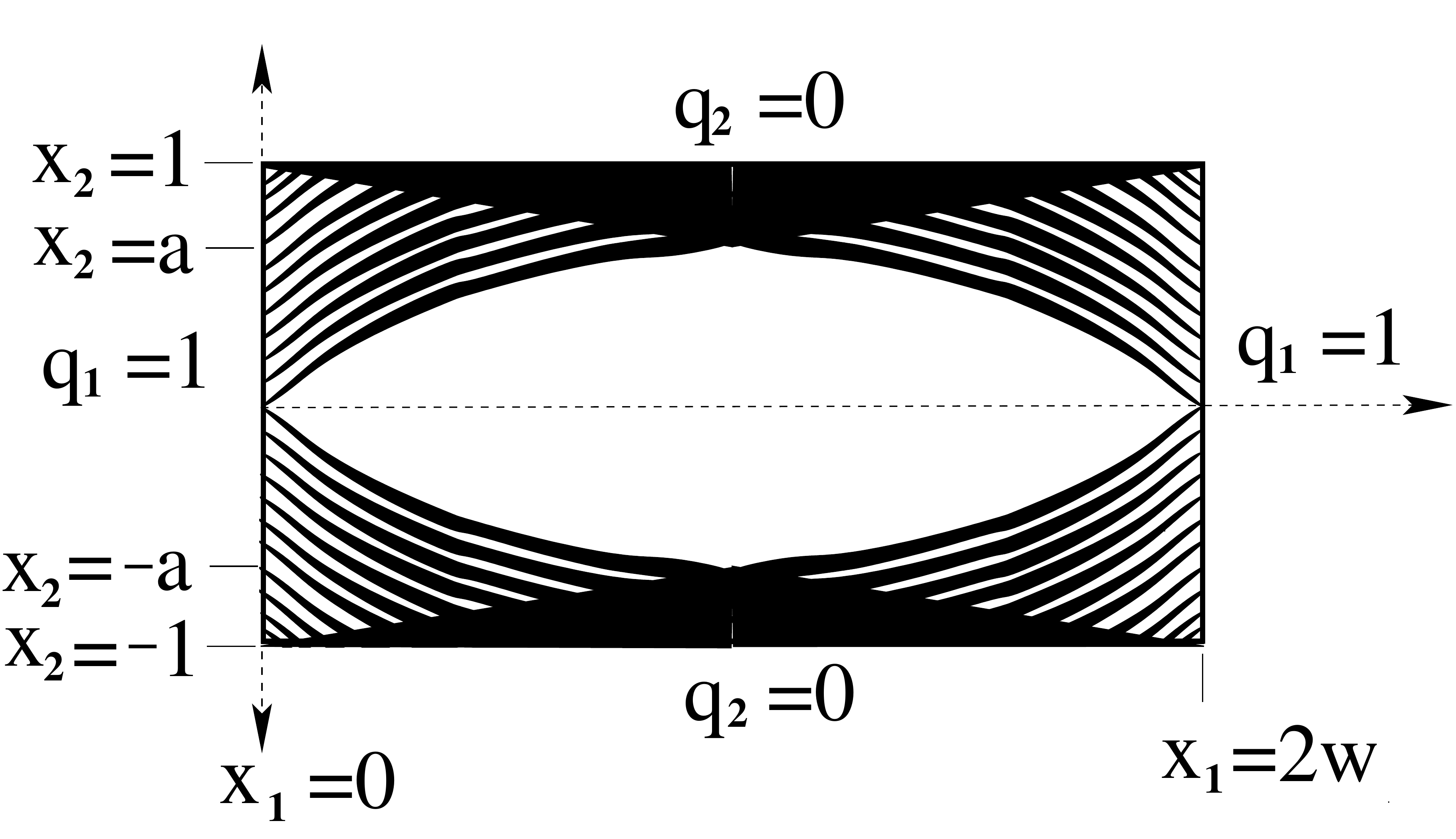}
\caption{A thermal shielding problem: distribute the two phases to minimize the heat flux through the window $x_1=w$ and $x_2\in[-a,a]$. Here the highly
conducting phase is the black one, while the poorer one is in white. The phase distribution shown is not the optimal one: the ribs need to be taken
thiner and thiner, approximating a laminate on the microscale. The optimal solution is provided in \protect\cite{Gibiansky:1988:OFH}.}
\labfig{2}
\end{figure}

\ajj{Now consider the special case where phase 2 is a thermally insulating material since, as remarked above (compare \eq{0.3a} and \eq{0.3c}),
it is most similar to the case of guiding stress using pentamodes.} The shielding problem of \fig{2} is then easily solved: just make sure that thermally insulating
material surrounds the window $x_1=1$ and $x_2\in[-a,a]$. To remove some of this degeneracy in the structure
one may require that phase 1 occupy a fixed volume fraction $p$ in the rectangle,
and additionally that the net thermal resistance,
\beq R=\int_{-1}^{+1}T(0,x_2)\,dx_2-\int_{-1}^{+1}T(2w,x_2)\,dx_2=\int_{0}^{2w}\,dx_1\int_{-1}^{+1}\Bq(x_1,x_2)\cdot\Be(x_1,x_2)\,dx_2,
\eeq{0.4a}
be minimized. Then at the microscale level one should use a laminate so that \eq{0.4} is attained for all $\Bx$ when $(\Bq_0,\Be_0)$ is replaced by 
$(\Bq(\Bx),\Be(\Bx))$ and $f$ is replaced by the local volume fraction $f(\Bx)$.
The task is then to consider all divergence free current fields $\Bq(\Bx)$, $\Div\Bq=0$, defined in the quarter rectangle $0\leq x_1\leq w$,
$0\leq x_2 \leq 1$, satisfying the boundary conditions
\beq q_1(0,x_2)=1\text{ for }0\leq x_2 \leq 1, \quad q_1(w,x_2)=0\text{ for } 0\leq x_2 \leq a, \quad q_2(x_1,0)=q_2(x_1,1)=0\text{ for }0\leq x_1 \leq w.
\eeq{0.4c}
Among these fields, and among all choices of $f(\Bx)$ with $1\geq f(\Bx)\geq 0$ that satisfy 
\beq p=\frac{1}{w}\int_{0}^{w}\,dx_1\int_{0}^{+1}f(x_1,x_2)\,dx_2, 
\eeq{0.4d}
the objective is to minimize
\beq  R=4\int_{0}^{w}\,dx_1\int_{0}^{+1}(f(\Bx)k_1)^{-1}|\Bq(\Bx)|^2\,dx_2.
\eeq{0.4e}
We use a Lagrange multiplier $\Gl^2$, with $\Gl>0$, for the constraint \eq{0.4e}. Then the optimal choice of $f(\Bx)$ is
\beqa f(\Bx)& = & \Gl|\Bq(\Bx)|\text{ if }|\Bq(\Bx)|<1/\Gl, \nonum
            & = & 1 \text{ otherwise},
\eeqa{0.4f}
where $\Gl$ depends on $\Bq$ and is chosen so \eq{0.4d} is satisfied.
Finally one optimizes over all $\Bq(\Bx)$, satisfying $\Div\Bq=0$ and \eq{0.4c}, the quantity
\beq R=\frac{4}{k_1\Gl}\int_{\GT}|\Bq(\Bx)|\,d\Bx+\frac{4}{k_1}\int_{\GT_C}|\Bq(\Bx)|^2, \eeq{0.4g}
where $\GT$ is the region of the quarter rectangle where $|\Bq(\Bx)|<1/\Gl$ and $\GT_C$
is the remaining region of the quarter rectangle where $|\Bq(\Bx)|\geq 1/\Gl$, and to satisfy \eq{0.4d},
\beq \Gl=(wp-|\GT_C|)/\int_{\GT}|\Bq(\Bx)|\,d\Bx, \eeq{0.4h}
in which $|\GT_C|$ is the area occupied by $\GT_C$. Assuming there
is a minimizer $\Bq(\Bx)$ one finds the associated temperature field $T(\Bx)$ by integrating 
\beq \Bq(x_1,x_2)\cdot\Be(x_1,x_2)=(f(\Bx)k_1)^{-1}|\Bq(\Bx)|^2 \eeq{0.4i}
back along the streamlines of $\Bq(\Bx)$, starting from the value $T(\Bx)=0$ when $x_1=w$. 

Of course similar sorts of optimization questions can be asked in the context of $3d$-elasticity rather than thermal conductivity.
Our complete determination of the weak $G$-closure for $3d$-printed materials opens the door for solving such problems and
an additional wide variety of important elasticity optimization problems for  $3d$-printed structures.
For instance, suppose we have a three-dimensional body $\GO$ containing a $3d$-printed structure made from a single isotropic material, labeled as phase 1, 
and that we are free to design its structure possibly with a constraint on its overall weight, or equivalently on the overall volume fraction $p$ occupied by phase 1.
Knowing the weak $G$-closure is sufficient if one wishes to optimize the response of the body in a single experiment where, for example, one prescribes the traction
$\Bt_0(\Bx)=\BGs(\Bx)\cdot\Bn$ at the surface $\Md\GO$ of $\GO$, where $\Bn$ is its outward normal, and optimizes the displacement 
$\Bu(\Bx)=\Bu_0(\Bx)$ at $\Md\GO$. Specifically, given a desired displacement $\Bu'_0(\Bx)$ and one may try to minimize the surface integral 
of $w(\Bx)|\Bu_0(\Bx)-\Bu'_0(\Bx)|^2$, where $w(\Bx)$ is some non-negative weighting function. (In a given geometry we are free to set $\Bu_0(\Bx)=\Bu'_0(\Bx)$ 
over any portion of the boundary where the void region intersects $\Md\GO$, but this can only be done on those portions of the boundary where 
the prescribed traction $\Bt_0(\Bx)$ is zero as void will not support a traction.)
Additionally one may want the structure to minimize the total compliance energy which can be identified with the integral over $\Md\GO$ of
$\Bu_0(\Bx)\cdot\Bt_0(\Bx)$. Allowing for body forces $\Bh(\Bx)$ the optimization problem for the given $\Bt_0(\Bx)$ becomes
\beq \inf_{\Gc(\Bx)}\int_{\Md\GO}w(\Bx)|\Bu_0(\Bx)-\Bu'_0(\Bx)|^2+k\Bu_0(\Bx)\cdot\Bt_0(\Bx)\,dS,
\eeq{0.5}
where $k>0$ and the minimization is over all geometries inside $\GO$, each associated with a characteristic functions $\Gc(\Bx)$, taking the value 1 
in phase 1, and zero otherwise, satisfying
\beq \int_\GO\Gc(\Bx)\,d\Bx=p, \eeq{0.5a}
and $\Bu_0(\Bx)$ is found by solving the elasticity equations
\beq \BGs(\Bx)=\BC_1\Gc(\Bx)\BGe(\Bx),\quad \Div\BGs=\Bh,\quad \BGe=[\Grad\Bu+(\Grad\Bu)^T]/2,
\eeq{0.6}
with the boundary condition $\BGs(\Bx)\cdot\Bn=\Bt_0(\Bx)$, for $\Bx\in\Md\GO$, and identifying $\Bu_0(\Bx)$ with the boundary
value of $\Bu(\Bx)$: $\Bu_0(\Bx)=\Bu(\Bx)$, for $\Bx\in\Md\GO$. When directly solving this optimization problem there might not be a minimizer
but rather a minimizing sequence with $\Gc(\Bx)$ having very rapid oscillations in certain regions of $\GO$. This would be indicative of the need to use composites in these regions
with microstructure much smaller than the dimensions of $\GO$. Then one can replace these composites with an homogenized material with elasticity tensor
matching that of the effective tensor. So to avoid oscillations, we follow the standard prescription in optimal design problems
(as described, for example, in the books \cite{Cherkaev:2000:VMS, Allaire:2002:SOH, Milton:2002:TOC, Bendsoe:2004:RPM, Tartar:2009:GTH})  
and reformulate the problem in terms of the relaxed problem:
\beq \inf_{\Bu(\Bx),\BGs(\Bx), f(\Bx)}\int_{\Md\GO}w(\Bx)|\Bu_0(\Bx)-\Bu'_0(\Bx)|^2+k\int_{\GO}W_{f(\Bx)}(\BGs(\Bx))\,d\Bx,
\eeq{0.7}
where the local volume fraction $f(\Bx)$ satisfies the constraints that $0\leq f(\Bx) \leq 1$ for all $\Bx\in\GO$, and
\beq  \int_\GO f(\Bx)\,d\Bx=p, \eeq{0.8}
while $\BGs(\Bx)$ satisfies the boundary conditions that $\BGs(\Bx)\cdot\Bn=\Bt_0(\Bx)$ for $\Bx\in\Md\GO$,
$\Bu_0(\Bx)$ is identified with the boundary value of $\Bu(\Bx)$: $\Bu_0(\Bx)=\Bu(\Bx)$, for $\Bx\in\Md\GO$,
and inside $\GO$ the fields satisfy the constraints
\beq \BGe(\Bx)\in\CE_{f(\Bx)}(\BGs(\Bx)),\quad \Div\BGs=\Bh,\quad \BGe=[\Grad\Bu+(\Grad\Bu)^T]/2.
\eeq{0.9}
\ajj{In other words, one replaces the constitutive law by the constraint that the pair $(\BGs(\Bx),\BGe(\Bx))$ lies in the weak $G$-closure associated
with the local volume fraction $f(\Bx)$, which implies there exists a local microstructure having $f(\Bx)$ as its volume fraction and with effective 
elasticity tensor $\BC_*(\Bx)$ such that $\BGs(\Bx)=\BC_*(\Bx)\BGe(\Bx)$. Due to the non-linearity of the constraint that $\BGe(\Bx)\in\CE_{f(\Bx)}(\BGs(\Bx))$ it is
evident that such problems are best solved numerically, and this is beyond the scope of the current paper.} 

Once the relaxed problem is solved one should, in principle, insert at each point $\Bx$ the composite that realizes the required (stress, strain) pair
$(\BGs(\Bx),\BGe(\Bx))$ at each point $\Bx\in\GO$. In practice one would look for a realistic single-scale microstructure that approximately realizes the 
pair $(\BGs(\Bx),\BGe(\Bx))$. We are not advocating the use of the complicated geometries described in \cite{Milton:2016:PEE}. Rather
we view the existence of complicated geometries that have the desired response as motivation for driving the numerical search for more
realistic single-scale microstructure that approximately achieve the same response.
This microstructure will change as $\Bx$ changes, but of course the change should be made in some sort of continuous
way that preserves the connectivity of phase 1. Other considerations may also weigh in. In particular, one may want the structure to not 
completely collapse when the boundary traction $\Bt(\Bx)$ is slightly perturbed. Accordingly, modifications of the microstructure might be needed
to take into account these additional requirements.

There are a wealth of other optimization problems that could be numerically solved using the weak $G$-closure. Instead
of prescribing the surface traction $\Bt_0(\Bx)$, one could fix the boundary displacement $\Bu_0(\Bx)$ and given
a desired traction $\Bt'_0(\Bx)$ try to minimize the surface integral of $w(\Bx)|\Bt_0(\Bx)-\Bt'_0(\Bx)|^2$ for some weight function $w(\Bx)\geq 0$.
Mixed boundary conditions are possible too. The characteristic feature is that we are optimizing the response in a single experiment, not for a series of
experiments where a given structure is exposed to multiple applied loadings.

%%%%%%%%%%%%%%%%%%%%%%%%%%%%%%%%%%%%%%%%%%%%%%%%%%%%%%%%%%%%%%%%%%%%%%%%%%%%%%%%%%%%%
\section{Establishing the main result}
\setcounter{equation}{0}
%%%%%%%%%%%%%%%%%%%%%%%%%%%%%%%%%%%%%%%%%%%%%%%%%%%%%%%%%%%%%%%%%%%%%%%%%%%%%%%%%%%%%%
The range of possible effective elasticity tensors that could be produced by mixing a given material having elasticity tensor $\BC_1$ with an extremely compliant phase
was explored in \cite{Milton:2016:PEE} where the following theorem \ajj{[in which we will only need the result pertaining to
$W_f^5(\BGs^0_1,\BGe^0_1,\BGe^0_2,\BGe^0_3,\BGe^0_4,\BGe^0_5)$]} was established:

\begin{theorem}(Milton, Briane and Harutyunyan)

Consider composites in three dimensions of two materials with positive definite elasticity tensors $\BC_1$ and $\BC_2=\Gd\BC_0$ mixed in proportions $f$ and $1-f$. Let the seven energy functions $W_f^k$, $k=0,1,\ldots,6$, that characterize the set $G_fU$ (with $U=(\BC_1,\Gd\BC_0)$) of possible elastic tensors $\BC^*$
be defined by
\beqa W_f^0(\BGs^0_1,\BGs^0_2,\BGs^0_3,\BGs^0_4,\BGs^0_5,\BGs^0_6)
& = &\min_{\BC^*\in GU_f}\sum_{j=1}^6\BGs^0_j:(\BC^*)^{-1}\BGs^0_j, \nonum
 W_f^1(\BGs^0_1,\BGs^0_2,\BGs^0_3,\BGs^0_4,\BGs^0_5,\BGe^0_1)
& = &\min_{\BC^*\in GU_f}\left[\BGe^0_1:\BC^*\BGe^0_1+\sum_{j=1}^5\BGs^0_j:(\BC^*)^{-1}\BGs^0_j\right],\nonum
 W_f^2(\BGs^0_1,\BGs^0_2,\BGs^0_3,\BGs^0_4,\BGe^0_1,\BGe^0_2)
& = &\min_{\BC^*\in GU_f}\left[\sum_{i=1}^2\BGe^0_i:\BC^*\BGe^0_i+\sum_{j=1}^4\BGs^0_j:(\BC^*)^{-1}\BGs^0_j\right],\nonum
 W_f^3(\BGs^0_1,\BGs^0_2,\BGs^0_3,\BGe^0_1,\BGe^0_2,\BGe^0_3)
& = &\min_{\BC^*\in GU_f}\left[\sum_{i=1}^3\BGe^0_i:\BC^*\BGe^0_i+\sum_{j=1}^3\BGs^0_j:(\BC^*)^{-1}\BGs^0_j\right],\nonum
 W_f^4(\BGs^0_1,\BGs^0_2,\BGe^0_1,\BGe^0_2,\BGe^0_3,\BGe^0_4)
& = &\min_{\BC^*\in GU_f}\left[\sum_{i=1}^4\BGe^0_i:\BC^*\BGe^0_i+\sum_{j=1}^2\BGs^0_j:(\BC^*)^{-1}\BGs^0_j\right],\nonum
 W_f^5(\BGs^0_1,\BGe^0_1,\BGe^0_2,\BGe^0_3,\BGe^0_4,\BGe^0_5)
& = & \min_{\BC^*\in GU_f}\left[\left(\sum_{i=1}^5\BGe^0_i:\BC^*\BGe^0_i\right)+\BGs^0_1:(\BC^*)^{-1}\BGs^0_1 \right],\nonum
 W_f^6(\BGe^0_1,\BGe^0_2,\BGe^0_3,\BGe^0_4,\BGe^0_5,\BGe^0_6)
& = & \min_{\BC^*\in GU_f}\sum_{i=1}^6\BGe^0_i:\BC^*\BGe^0_i.
\eeqa{1.1}
\ajj{Each energy function here represents the sum of six elastic energies, each obtained from an experiment where the
composite, with effective tensor $\BC^*$, is either subject to an applied stress $\BGs^0_i$ or an applied strain $\BGe^0_j$. 
These 7 functions characterize the $G$-closure in much the same way that a convex set is characterized by its Legendre transform.}
The applied strains $\BGe^0_i$ and applied stresses $\BGs^0_j$ can be chosen to meet the orthogonality condition 
\beq \BGe^0_i:\BGs^0_j  =  0,\quad \BGe^0_i:\BGe^0_k=0,\quad \BGs^0_j:\BGs^0_\ell=0 \quad{\rm for~all~}i, j, k,\ell~{\rm with}~i\ne j, ~
i\ne k, ~j\ne \ell.
\eeq{1.2}
The energy function $W_f^0$ is given by
\beq W_f^0( \BGs^0_1, \BGs^0_2, \BGs^0_3, \BGs^0_4, \BGs^0_5, \BGs^0_6)=\sum_{j=1}^6\BGs^0_j:\BC_f^A(\BGs^0_1,\BGs^0_2,\BGs^0_3,\BGs^0_4,\BGs^0_5,\BGs^0_6)\BGs^0_j,
\eeq{1.3}
as proved by Avellaneda \cite{Avellaneda:1987:OBM}. Here
$\BC_f^A(\BGs^0_1,\BGs^0_2,\BGs^0_3,\BGs^0_4,\BGs^0_5,\BGs^0_6)$ is the effective elasticity tensor of a complementary Avellaneda material, that is a sequentially layered laminate with the minimum value of the sum of complementary energies 
\beq \sum_{j=1}^6\BGs^0_j:(\BC^*)^{-1}\BGs^0_j.=\BA\Shortstack{. . . .}(\BC^*)^{-1},\quad \BA=\sum_{j=1}^6\BGs^0_j\otimes\BGs^0_j,
\eeq{1.4}
where the four vertical dots represent the contraction of the four indices of $(\BC^*)^{-1}$ with the four indices of $\BA$.
Additionally, \ajj{by introducing into the homogenized Avellaneda material appropriate sets of thin walls, in different directions, that themselves
are microstructured in such a way that they support some loadings, but slip or easily compress or easily expand under other loadings,}
we now have
\beqa
\lim_{\Gd\to 0}W_f^3(\BGs^0_1,\BGs^0_2,\BGs^0_3,\BGe^0_1,\BGe^0_2,\BGe^0_3) & = & \sum_{j=1}^3\BGs^0_j:[\BC_f^A(\BGs^0_1,\BGs^0_2,\BGs^0_3,0,0,0)]^{-1}\BGs^0_j ,\nonum
\lim_{\Gd\to 0}W_f^4(\BGs^0_1,\BGs^0_2,\BGe^0_1,\BGe^0_2,\BGe^0_3,\BGe^0_4) & = & \sum_{j=1}^2\BGs^0_j:[\BC_f^A(\BGs^0_1,\BGs^0_2,0,0,0,0)]^{-1}\BGs^0_j , \nonum
\lim_{\Gd\to 0}W_f^5(\BGs^0_1,\BGe^0_1,\BGe^0_2,\BGe^0_3,\BGe^0_4,\BGe^0_5) & = &  \BGs^0_1:[\BC_f^A(\BGs^0_1,0,0,0,0,0)]^{-1}\BGs^0_1 , \nonum
\lim_{\Gd\to 0}W_f^6(\BGe^0_1,\BGe^0_2,\BGe^0_3,\BGe^0_4,\BGe^0_5,\BGe^0_6) & = & 0,
\eeqa{1.5}
for all combinations of applied stresses $\BGs_j^0$ and applied strains $\BGe^0_i$ meeting the orthogonality conditions \eq{1.2}.
\ajj{Thus the material in the walls that are inserted into the homogenized Avellaneda material
is microstructured in such a way that the applied strains $\BGe^0_j$ in these energy functions now cost no energy, while 
there is negligible effect on the energy associated with the applied stresses $\BGs^0_i$, as the walls support these loadings.}
When  $\det\BGe^0_1=0$ but $\BGe^0_1$ is not positive or negative semidefinite, we have
\beq
W_f^1(\BGs^0_1,\BGs^0_2,\BGs^0_3,\BGs^0_4,\BGs^0_5,\BGe^0_1)=\sum_{j=1}^5\BGs^0_j:[\BC_f^A(\BGs^0_1,\BGs^0_2,\BGs^0_3,\BGs^0_4,\BGs^0_5,0)]^{-1}\BGs^0_j,
\eeq{1.6}
while when the equation $\det(\BGe^0_1+t\BGe^0_2)=0$ has at least two distinct roots for $t$, and additionally, the matrix pencil $\BGe(t)=\BGe^0_1+t\BGe^0_2$ does not contain any positive definite or negative definite matrices as $t$ is varied, we have 
\beq 
W_f^2(\BGs^0_1,\BGs^0_2,\BGs^0_3,\BGs^0_4,\BGe^0_1,\BGe^0_2)=\sum_{j=1}^4\BGs^0_j:[\BC_f^A(\BGs^0_1,\BGs^0_2,\BGs^0_3,\BGs^0_4,0,0)]^{-1}\BGs^0_j.
\eeq{1.7}
\end{theorem}

As in \cite{Milton:2016:PEE} we remark that the material with tensor $\BC_2$ could itself be formed from the material with tensor $\BC_1$ with sufficiently many disconnected voids carved out
so that its effective elasticity tensor, that we equate with $\BC_2$, is extremely small. Of course its microstructure must be much smaller than the microstructure in the composites of
$\BC_1$ and $\BC_2$ that are considered. \ajj{ Also we emphasize that this theorem is only valid for linear elasticity with infinitesimal displacements. Indeed the microstructures required 
are quite delicate having structure on multiple length scales. Approximations to these structures may function for small vibrations, but are unlikely to function well under
large deformations.}

\ajj{An additional remark is that ``the complementary Avellaneda material'' is not uniquely determined. Generally it is constructed by taking a laminate of both phases, and then sequentially
laminating this on larger and larger length scales with the stiffer phase in a variety of different directions (the consideration of the effective properties of such
multiscale laminates goes back to Maxwell \cite{Maxwell:1954:TEMb}). An explicit formula for the resulting effective elasticity tensor, in terms of 
the lamination directions and volume fractions of the stiffer phase introduced at each level has been given by Francfort and Murat \cite{Francfort:1986:HOB}. From their formulae
it is evident that the ordering in which this layering occurs does not matter and that many such hierarchical laminates produce the same effective tensor.
To obtain  ``the complementary Avellaneda material'' one needs to numerically optimize the lamination directions and volume fractions of the stiffer phase introduced at each level.
In the case when $\BA$ in \eq{1.4} is an isotropic fourth order tensor, ``the complementary Avellaneda material'' can be taken as any elastically isotropic material 
\cite{Norris:1985:DSE,Francfort:1986:HOB,Milton:1986:MPC}
that simultaneously attains the upper Hashin-Shtrikman bulk and shear modulus bounds \cite{Hashin:1963:VAT}.
}

\ajj{We also point out that the microstructured material within the walls (that are inserted into the Avellaneda material) can be chosen so it has precisely the
desired volume fraction $f$, if necessary by adding islands of phase 1 (possibly with thin ligaments attaching them to the main structure). Each
set of parallel walls are easily strained under a strain or strains that are necessarily of the form $(\Ba\otimes\Bn+\Bn\otimes\Ba)/2$, where $\Bn$
is the normal to the walls. The applied strains entering \eq{1.5}, \eq{1.6} and \eq{1.7} that are associated with $W_f^i$ must each be a linear superposition 
of $i$ such symmetrized rank one tensors. Thus the crux of the theorem is identifying those $i$-dimensional subspaces in the $6$-dimensional space
of symmetric $3\times 3$ matrices that can be spanned by $i$ symmetrized rank-one matrices. The equalities \eq{1.5} reflect the fact that any   
$i$-dimensional subspace can be spanned in this way if $i\geq 3$. The inequalities \eq{1.6} and \eq{1.7} reflect the fact that some restrictions on the subspace
are needed if $i<3$. In the case $i=1$, $\BGe^0_1$ must itself be a symmetrized rank one matrix and this is equivalent to the conditions
that  $\det\BGe^0_1=0$ but $\BGe^0_1$ is not positive or negative semidefinite.}

In the analysis in this paper we assume $\BC_1$ is an isotropic elasticity tensor, defined by its action
\beq \BC_1\BGe=2\Gm\BGe+\Gl(\Tr\BGe)\BI,
\eeq{1.7a}
where $\BI$ represents the $3\times 3$ identity matrix, $\Gm$ is the shear modulus and $\Gl$ is the Lame modulus.
As we will only be interested in the function $W_f^5(\BGs^0_1,\BGe^0_1,\BGe^0_2,\BGe^0_3,\BGe^0_4,\BGe^0_5)$, the relevant part of the theorem implies
\beq 
\lim_{\Gd\to 0}W_f^5(\BGs^0,\BGe^0_1,\BGe^0_2,\BGe^0_3,\BGe^0_4,\BGe^0_5)  =   \BGs^0:[\BC_f^{A}(\BGs^0)]^{-1}\BGs^0\equiv W_f(\BGs^0),
\eeq{1.8}
where $\BC_f^A(\BGs^0)$ is the effective elasticity tensor of a complementary Avellaneda material, that is a sequentially layered laminate with the minimum value of the complementary energy $\BGs^0:(\BC^*)^{-1}\BGs^0$ among all composites. An explicit expression for the function $W_f(\BGs^0)$ in the usual case where $\Gl>0$ was obtained independently by
Cherkaev and Gibiansky \cite{Gibiansky:1987:MCE} and Allaire \cite{Allaire:1994:ELP}, and (using mostly Allaire's notation) is given by 
\beq W_f(\BGs^0)=\BGs^0:\BC_1^{-1}\BGs^0+\frac{f}{2\Gm}g(\BC_1,\BGs^0), \eeq{1.10}
where $g(\BC_1,\BGs^0)$ depends on the Lame moduli, $\Gl$ and $\Gm$ of $\BC_1$ and on the eigenvalues $\Gs_1$, $\Gs_2$, and $\Gs_3$ of $\BGs^0$
that we can assume to be labeled so $\Gs_1\leq\Gs_2\leq\Gs_3$. As $W_f(\BGs^0)=W_f(-\BGs^0)$ it suffices to assume that $\BGs^0$ has at most one negative 
eigenvalue, and in the case all eigenvalues are nonnegative
\beqa
g(\BC,\BGs) & = & \frac{2\Gm+\Gl}{2(2\Gm+3\Gl)}(\Gs_1+\Gs_2+\Gs_3)^2\text{ if }\Gs_3\leq \Gs_1+\Gs_2, \nonum
            & = & (\Gs_1+\Gs_2)^2+\Gs_3^2-\frac{\Gl}{2\Gm+3\Gl}(\Gs_1+\Gs_2+\Gs_3)^2\text{ if }\Gs_3\geq \Gs_1+\Gs_2,
\eeqa{1.12}
while when one eigenvalue, namely $\Gs_1$, is negative,
\beqa
g(\BC,\BGs) & = & \frac{2\Gm+\Gl}{2(2\Gm+3\Gl)}\left(\Gs_3+\Gs_2-\frac{\Gm+2\Gl}{\Gm+\Gl}\Gs_1\right)^2 \nonum
& ~& \text{ if } \Gs_3+\Gs_2\geq \frac{-\Gm}{\Gm+\Gl}\Gs_1\text{ and } \Gs_3-\Gs_2\leq \frac{-\Gm}{\Gm+\Gl}\Gs_1, \nonum
 & = & (\Gs_3+\Gs_2)^2+\Gs_1^2-\frac{\Gl}{2\Gm+3\Gl}(\Gs_1+\Gs_2+\Gs_3)^2 \text{ if } \Gs_3+\Gs_2\leq \frac{-\Gm}{\Gm+\Gl}\Gs_1, \nonum
& = & \Gs_1^2+\Gs_2^2+\Gs_3^2-\frac{2\Gm}{\Gm+\Gl}\Gs_1\Gs_2-\frac{\Gl}{2\Gm+3\Gl}(\Gs_1+\Gs_2+\Gs_3)^2
\text{ if }  \Gs_3-\Gs_2\geq \frac{-\Gm}{\Gm+\Gl}\Gs_1.
\eeqa{1.15}
Given any $\BGs^0$, Cherkaev, Gibiansky, and Allaire also provide the precise prescription for building \ajj{the relevant Avellaneda material} namely
a sequential laminate of rank 2 or 3 that attains the bounds.
Rank 1 means layers of phase 1 separated by void. The needed rank 2 materials are obtained by slicing this rank 1 laminate into layers, orders of magnitude thicker than the layers
in the rank 1 material and with a different orientation of the cut, that are then separated by  
layers of phase 1. The needed rank 3 materials are obtained by slicing this rank 2 laminate into layers, orders of magnitude thicker than the layers
in the rank 2 material, that are then separated again by layers of phase 1.

\ajj{Our goal is to show that for any given $\BGs^0$ and $\BGe^0$ satisfying \eq{0.1} as an inequality, there exists a composite of the two materials with
tensors $\BC_1$ and $\BC_2=\Gd\BC_0$ mixed in proportions $f$ and $1-f$ such that if $\BGe^0$ is the average strain in the material, then the 
average stress in the material can be made arbitrarily close to $\BGs^0$ for a sufficiently small value of $\Gd$. Heuristically, suppose $\Gd=0$ and 
we have an optimal pentamode with elasticity tensor
\beq \BC^*_0=\Gg\BGs^0\otimes\BGs^0, \eeq{1.15a}
where by optimal we mean that for an applied stress $\BGs^0$ there exists a strain $\widetilde{\BGe}^0$ such that $\BGs^0=\BC^*_0\widetilde{\BGe}^0$
and one has equality in \eq{0.1} (with $\BGe^0$ replaced by $\widetilde{\BGe}^0$) , i.e. $\Gg=W_f(\BGs^0)/(\BGs_0:\BGs_0)$. Then clearly $\BC^*_0\BGe^0$ is
our applied stress $\BGs^0$ and so $\BGe^0$ can be identified with an applied strain that generates the stress $\BGs^0$ (with $\BGe^0-\widetilde{\BGe}^0$ lying in the null space
of $\BC^*_0$). The argument seems simple, but complications arise if we take $\Gd$ non-zero and take the limit $\Gd\to 0$. Then, given an
applied stress $\BGs^0$ the associated average strain $\widetilde{\BGe}^0=(\BC^*_\Gd)^{-1}\BGs^0$ is uniquely determined and not necessarily close
to the desired average strain $\BGe^0$. The way around this difficulty is to prescribe the average strain to be $\BGe^0$ and then show that the average
stress converges to $\BGs^0$. We remark, however, that while the bound \eq{0.1} is a linear constraint on $\BGe^0$ given $\BGs^0$, it is a highly
non-trivial, non-linear constraint on $\BGs^0$ given $\BGe^0$. Obviously, care is required and to make things rigorous we need to use a more sophisticated
argument, which we now proceed in doing.}  

Let us set $\BGe^0_0=\BGs^0/t$ and let us choose $t$ and the matrices $\BGe^0_1$, $\BGe^0_2$, $\BGe^0_3$, $\BGe^0_4$, $\BGe^0_5$ so that together with $\BGe^0_0$ they form
an orthonormal basis on the space of $3\times 3$ symmetric matrices. Then \eq{1.8} implies there exists a sequence $c_\Gd$ of positive numbers
with $c_\Gd\to 0$ as $\Gd\to 0$, and a sequence of positive definite effective elasticity tensors $\BC^*_\Gd$  associated with a sequence of
two-phase microgeometries (made precise in \cite{Milton:2016:PEE}) with the phases having tensors $\BC_1$ and $\BC_2=\Gd\BC_0$,
such that
\beq W_f(\BGs^0)\leq \left(\sum_{i=1}^5\BGe^0_i:\BC^*_\Gd\BGe^0_i\right)+\BGs^0:(\BC^*_\Gd)^{-1}\BGs^0 \leq W_f(\BGs^0)+c_\Gd.
\eeq{1.16}
The elastic complementary energy bound \eq{0.1} then implies
\beq 0\leq \BGs^0:(\BC^*_\Gd)^{-1}\BGs^0-W_f(\BGs^0)\leq c_\Gd,\quad 0\leq\sum_{i=1}^5\BGe^0_i:\BC^*_\Gd\BGe^0_i\leq c_\Gd.
\eeq{1.17}
Using the basis $\BGe^0_0$, $\BGe^0_1$, $\BGe^0_2$, $\BGe^0_3$, $\BGe^0_4$, $\BGe^0_5$ the elasticity tensor $\BC^*_\Gd$ takes the form
\beq \BC^*_\Gd=\bpm \Ga_\Gd & \Ba_\Gd^T \\
                   \Ba_\Gd & \BA_\Gd \epm,
\eeq{1.18}
where for each $\Gd>0$, $\Ga_\Gd$ is a positive scalar, $\Ba_\Gd$ is a 5-dimensional vector, and $\BA_\Gd$ is a positive definite $5\times 5$ matrix. The matrix \eq{1.18} has inverse
\beq (\BC^*_\Gd)^{-1}=\bpm (\Ga_\Gd-\Ba_\Gd\cdot\BA_\Gd^{-1}\Ba_\Gd)^{-1} & -(\Ga_\Gd-\Ba_\Gd\cdot\BA_\Gd^{-1}\Ba_\Gd)^{-1}\Ba_\Gd^T\BA_\Gd^{-1} \\
-\BA_\Gd^{-1}\Ba_\Gd(\Ga_\Gd-\Ba_\Gd\cdot\BA_\Gd^{-1}\Ba_\Gd)^{-1} & \BA_\Gd^{-1}+\frac{\BA_\Gd^{-1}\Ba_\Gd\Ba_\Gd^T\BA_\Gd^{-1}}{(\Ga_\Gd-\Ba_\Gd\cdot\BA_\Gd^{-1}\Ba_\Gd)},
\epm
\eeq{1.23}
and so the first inequality in \eq{1.17}, together with the positivity of $\BA_\Gd$, implies
\beq  W_f(\BGs^0)+c_\Gd \geq t^2(\Ga_\Gd-\Ba_\Gd\cdot\BA_\Gd^{-1}\Ba_\Gd)^{-1}\geq t^2/\Ga_\Gd.
\eeq{1.24}
Hence $1/\Ga_\Gd$ remains bounded as $\Gd\to 0$.
The second bound in \eq{1.17} and the orthonormality of the matrices $\BGe^0_i$, $i=1,2,\ldots,5$, imply $\BA_\Gd\leq c_\Gd\BI$. As $\BC^*_\Gd$ is a positive definite matrix,
the Schur complement $\BA_\Gd-\Ga_\Gd^{-1}\Ba_\Gd\Ba_\Gd^T$ must also be positive definite and so we have
\beq 0\leq \Ba_\Gd\cdot(\BA_\Gd-\Ga_\Gd^{-1}\Ba_\Gd\Ba_\Gd^T)\Ba_\Gd \leq \Ba_\Gd\cdot(c_\Gd\BI-\Ga_\Gd^{-1}\Ba_\Gd\Ba_\Gd^T)\Ba_\Gd=c_\Gd|\Ba_\Gd|^2-\Ga_\Gd^{-1}|\Ba_\Gd|^4,
\eeq{1.25}
implying $c_\Gd>\Ga_\Gd^{-1}|\Ba_\Gd|^2$. So if $c_\Gd$ is small, so too is $\Ba_\Gd$ as $1/\Ga_\Gd$ remains bounded. Now \eq{1.18} implies
\beq \bpm  t \\  t\Ba_\Gd/\Ga_\Gd  \epm =\BC^*_\Gd\bpm t/\Ga_\Gd \\ 0 \epm, 
\eeq{1.25a}
giving
\beq t^2/\Ga_\Gd=\bpm  t \\  t\Ba_\Gd/\Ga_\Gd  \epm\cdot (\BC^*_\Gd)^{-1}\bpm  t \\  t\Ba_\Gd/\Ga_\Gd  \epm \geq W_f(\BGs^0+\BGs_R),
\eeq{1.26}
with
\beq \BGs_R=\sum_{i=1}^5 t\Ga_\Gd^{-1}\{\Ba_\Gd\}_i\BGe^0_i,
\eeq{1.27}
approaching zero as $\Gd\to 0$. With \eq{1.24} we deduce that 
\beq  W_f(\BGs^0)+c_\Gd \geq t^2/\Ga_\Gd \geq W_f(\BGs^0+\BGs_R).
\eeq{1.28}
As the energy function $W_f(\BGs^0)$ is a continuous function of $\BGs^0$, and since $\BGs_R\to 0$ as $\Gd\to 0$, we conclude that
\beq \Ga_\Gd\to t^2/W_f(\BGs^0)\text{ as }\Gd\to 0.
\eeq{1.29}
Now $\BGs^0$ and any given $\BGe^0$ such that $\BGe^0:\BGs^0=W_f(\BGs^0)$ are represented in our basis by vectors of the form
\beq \BGs^0=\bpm t \\ 0 \epm,\quad\BGe^0=\bpm W_f(\BGs^0)/t \\ \BGe^0_\perp \epm,
\eeq{1.30} 
where the 5-component vector $\BGe^0_\perp$ represents the component of $\BGe^0$ perpendicular to $\BGs^0$ represented in the basis $\BGe^0_i$, $i=1,2,\ldots,5$. Accordingly, using \eq{1.18},
$\BGs_\Gd^0=\BC^*_\Gd\BGe^0$ is represented in our basis by
\beq \BGs_\Gd^0=\bpm \Ga_\Gd W_f(\BGs^0)/t +\Ba_\Gd\cdot\BGe^0_\perp \\  \Ba_\Gd W_f(\BGs^0)/t +\BA_\Gd\BGe^0_\perp \epm.
\eeq{1.31}
As both $\Ba_\Gd\to 0$ and $0\leq \BA_\Gd\leq c_\Gd\BI\to 0$ as $\Gd\to 0$ and \eq{1.29} holds, we see that $\BGs_\Gd^0$ converges to $\BGs^0$ as $\Gd\to 0$. 

%%%%%%%%%%%%%%%%%%%%%%%%%%%%%%%%%%%%%%%%%%%%%%%%%%%%%%%%%%%%%%%%%%%%%%%%%%%%%%%%%%%%%%%%%%%%%%%%%%%%%%%%%%%%%%%%%%%%%%%%%%%%%%%%%
\section{Near optimal unimodes for guiding strain while minimizing stiffness in mixtures of a compliant and a rigid material}
\setcounter{equation}{0}
%%%%%%%%%%%%%%%%%%%%%%%%%%%%%%%%%%%%%%%%%%%%%%%%%%%%%%%%%%%%%%%%%%%%%%%%%%%%%%%%%%%%%%%%%%%%%%%%%%%%%%%%%%%%%%%%%%%%%%%%%%%%%%
The range of possible effective elasticity tensors that could be produced by mixing a given material having elasticity tensor $\BC_1$ with an almost rigid phase
was explored in \cite{Milton:2016:TCC} where the following theorem \ajj{[in which we will only need the result pertaining to
$W_f^1(\BGs^0_1,\BGs^0_2,\BGs^0_3,\BGs^0_4,\BGs^0_5,\BGe^0_1)$]} was established:

\begin{theorem}(Milton, Harutyunyan, and Briane)

Consider composites in three dimensions of two materials with positive definite elasticity tensors $\BC_1$ and $\BC_2=\Gd\BC_0$ mixed in proportions $f$ and $1-f$. Let the seven energy functions $W_f^k$, $k=0,1,\ldots,6$, that characterize the set $G_fU$ (with $U=(\BC_1,\Gd\BC_0)$) of possible elastic tensors $\BC^*$
be defined by \eq{1.1}. These  energy functions involve a set of applied strains $\BGe^0_i$ and applied stresses $\BGs^0_j$ meeting the orthogonality condition \eq{1.2}. The energy function $W_f^6$ is given by
\beq
 W_f^6(\BGe^0_1,\BGe^0_2,\BGe^0_3,\BGe^0_4,\BGe^0_5,\BGe^0_6)
 = \sum_{i=1}^6\BGe^0_i:\widetilde{\BC}_f^A(\BGe^0_1,\BGe^0_2,\BGe^0_3,\BGe^0_4,\BGe^0_5,\BGe^0_6)\BGe^0_i,
\eeq{2.1}
as established by Avellaneda \cite{Avellaneda:1987:OBM}, where $\widetilde{\BC}_f^A(\BGe^0_1,\BGe^0_2,\BGe^0_3,\BGe^0_4,\BGe^0_5,\BGe^0_6)$ is the effective elasticity tensor of an Avellaneda material, that is a sequentially layered laminate with the minimum value of the sum of elastic energies 
\beq \sum_{i=1}^6\BGe^0_j:\BC^*\BGe^0_j.
\eeq{2.2}
Again some of the applied stresses $\BGs_j^0$ or applied strains $\BGe^0_i$ could be zero. Additionally we now have
\beqa
 \lim_{\Gd\to \infty}W_f^0(\BGs^0_1,\BGs^0_2,\BGs^0_3,\BGs^0_4,\BGs^0_5,\BGs^0_6) & = & 0,\nonum
 \lim_{\Gd\to \infty}W_f^1(\BGs^0_1,\BGs^0_2,\BGs^0_3,\BGs^0_4,\BGs^0_5,\BGe^0_1) & = & \BGe^0_1:[\widetilde{\BC}_f^A(0, 0, 0, 0, 0, \BGe^0_1)]\BGe^0_1, \nonum
\lim_{\Gd\to \infty}W_f^2(\BGs^0_1,\BGs^0_2,\BGs^0_3,\BGs^0_4,\BGe^0_1,\BGe^0_2)   & = & \sum_{i=1}^2\BGe^0_i:[\widetilde{\BC}_f^A(0, 0, 0, 0,\BGe^0_1,\BGe^0_2)]\BGe^0_i, \nonum
\lim_{\Gd\to \infty} W_f^3(\BGs^0_1,\BGs^0_2,\BGs^0_3,\BGe^0_1,\BGe^0_2,\BGe^0_3)   &= &\sum_{i=1}^3\BGe^0_i:[\widetilde{\BC}_f^A(0, 0, 0, \BGe^0_1,\BGe^0_2,\BGe^0_3)]\BGe^0_i,
\eeqa{2.3}
for all combinations of applied stresses $\BGs_j^0$ and applied strains $\BGe^0_i$. When $\det(\BGs^0_1)=0$ we have
\beq
\lim_{\Gd\to \infty} W_f^5(\BGs^0_1,\BGe^0_1,\BGe^0_2,\BGe^0_3,\BGe^0_4,\BGe^0_5)=\sum_{i=1}^5\BGe^0_i:[\widetilde{\BC}_f^A(0,\BGe^0_1,\BGe^0_2,\BGe^0_3,\BGe^0_4,\BGe^0_5)]\BGe^0_i,
\eeq{2.4}
while, when $f(t)=\det(\BGs^0_1+t\BGs^0_2)=0$ has at least two roots, 
\beq
\lim_{\Gd\to \infty}W_f^4(\BGs^0_1,\BGs^0_2,\BGe^0_1,\BGe^0_2,\BGe^0_3,\BGe^0_4)=\sum_{i=1}^4\BGe^0_i:[\widetilde{\BC}_f^A(0, 0, \BGe^0_1,\BGe^0_2,\BGe^0_3,\BGe^0_4)]\BGe^0_i.
\eeq{2.5}
\end{theorem}
We remark that the material with tensor $\BC_2$ could itself be formed from the material with tensor $\BC_1$ with sufficiently many rigid inclusions inserted so 
that its effective elasticity tensor, that we equate with $\BC_2$,  is extremely large. Of course its microstructure must be much smaller than the microstructure in the composites of
$\BC_1$ and $\BC_2$ that are considered.

As we will only be interested in the function $W_f^1(\BGs^0_1,\BGs^0_2,\BGs^0_3,\BGs^0_4,\BGs^0_5,\BGe^0_1)$, the relevant part of the theorem implies
\beq \lim_{\Gd\to \infty}W_f^1(\BGs^0_1,\BGs^0_2,\BGs^0_3,\BGs^0_4,\BGs^0_5,\BGe^0)=\BGe^0:\widetilde{\BC}_f^{A}(\BGe^0)\BGe^0\equiv \widetilde{W}_f(\BGe^0),
\eeq{2.6}
where $\widetilde{\BC}_f^A(\BGe^0)$ is the effective elasticity tensor of an Avellaneda material, that is a sequentially layered laminate with the minimum value of the 
elastic energy $\BGe^0:\BC^*\BGe^0$ among all composites. We assume that the tensor $\BC_1$ is elastically isotropic with action given by
\eq{1.7a} with Lame parameter $\Gl>0$. In this case, Allaire and Kohn \cite{Allaire:1993:OBE} established that the function $\widetilde{W}_f(\BGe^0)$ is given by
\beq  \widetilde{W}_f(\BGe^0)= \BGe^0:\BC_1\BGe^0 +(1-f)\max_{\BGn}[2\BGe^0 :\BGn-f g(\BGn)], 
\eeq{2.7}
where $g(\BGn)$ is a function of the eigenvalues $\Gn_1, \Gn_2$,
and $\Gn_3$ of the symmetric matrix $\BGn$. Assuming that
these are labeled with $$ \Gn_1\leq\Gn_2\leq\Gn_3, $$
they provide the formula
\beqa
g(\BGn)&=&
\frac{(\Gn_1-\Gn_3)^2}{4\Gm}+\frac{(\Gn_1+\Gn_3)^2}{4(\Gl+\Gm)}  ~~{\rm if
}~~
\Gn_3\geq\frac{\Gl+2\Gm}{2(\Gl_2+\Gm_2)}(\Gn_1+\Gn_3)\geq\Gn_1, \nonum
g(\BGn)&=&
\frac{\Gn_1^2}{\Gl+2\Gm}  ~~{\rm if}~~
\Gn_1>\frac{\Gl+2\Gm}{2(\Gl+\Gm)}(\Gn_1+\Gn_3), \nonum
g(\BGn)&=&
\frac{\Gn_3^2}{\Gl+2\Gm}  ~~{\rm if}~~
\Gn_3<\frac{\Gl+2\Gm}{2(\Gl+\Gm)}(\Gn_1+\Gn_3).
\eeqa{2.8}
Let us set $\BGs^0_0=\BGe^0/t$ and let us choose $t$ and the matrices $\BGs^0_1$, $\BGs^0_2$, $\BGs^0_3$, $\BGs^0_4$, $\BGs^0_5$ so that together with $\BGs^0_0$ they form
an orthonormal basis on the space of $3\times 3$ symmetric matrices. Then it follows from \eq{2.6} that there exists a sequence $c_\Gd$ of positive numbers
with $c_\Gd\to 0$ as $\Gd\to \infty$, and a sequence of positive definite effective compliance tensors $\BS^*_\Gd$  associated with a sequence of
two-phase microgeometries (made precise in \cite{Milton:2016:PEE}) with the phases having compliance tensors $\BS_1=\BC_1^{-1}$ and $\BS_2=\BC_2^{-1}=\Gd^{-1}\BC_0^{-1}$,
such that
\beq \widetilde{W}_f(\BGe^0)\leq \left(\sum_{i=1}^5\BGs^0_i:\BS^*_\Gd\BGs^0_i\right)+\BGe^0:(\BS^*_\Gd)^{-1}\BGe^0 \leq \widetilde{W}_f(\BGe^0)+c_\Gd,
\eeq{2.9}
for all $\Gd>0$. The elastic energy bound then implies
\beq 0\leq \BGe^0:(\BS^*_\Gd)^{-1}\BGe^0-\widetilde{W}_f(\BGe^0)\leq c_\Gd,\quad 0\leq\sum_{i=1}^5\BGs^0_i:\BS^*_\Gd\BGs^0_i\leq c_\Gd.
\eeq{2.10}
Using the basis $\BGs^0_0$, $\BGs^0_1$, $\BGs^0_2$, $\BGs^0_3$, $\BGs^0_4$, $\BGs^0_5$ the compliance tensor $\BS^*_\Gd$ takes the form
\beq \BS^*_\Gd=\bpm \Gb_\Gd & \Bb_\Gd^T \\
                   \Bb_\Gd & \BB_\Gd \epm,
\eeq{2.11}
where $\Gb_\Gd$ is a scalar, $\Bb_\Gd$ is a 5-dimensional vector, and $\BB_\Gd$ is a $5\times 5$ matrix. The subsequent
analysis (which, for brevity, we omit) proceeds exactly in parallel to that in the previous section (with the roles of stresses and strains
interchanged, and with compliance tensors replacing elasticity tensors). We conclude that given any
$\BGs^0$ and $\BGe^0$ such that equality holds in the energy bounds \eq{0.10} it follows that 
$\BGe_\Gd^0=\BS^*_\Gd\BGs^0$ converges to $\BGe^0$ as $\Gd\to \infty$.

%%%%%%%%%%%%%%%%%%%%%%%%%%%%%%%%%%%%%%%%%%%%%%%%%%%%%%%%%%%%%%%%%%%%%%%%%%%%%%%%%%%
\section{Conclusion}
\setcounter{equation}{0}
%%%%%%%%%%%%%%%%%%%%%%%%%%%%%%%%%%%%%%%%%%%%%%%%%%%%%%%%%%%%%%%%%%%%%%%%%%%%%%%%%%%%%%
We have established the following two theorems:

\begin{theorem}
For $3d$-printed materials constructed from phase 1 with elasticity tensor $\BC_1$ defined by \eq{1.7a}
and void, with porosity $1-f$, the complementary energy bound
\beq W_f(\BGs^0)\leq\BGs^0:\BGe^0, \eeq{3.1}
where the explicit formula for $W_f(\BGs^0)$ is given by \eq{1.10}-\eq{1.15},
not only provides a sharp lower bound on the complementary energy, as established 
by Gibiansky and Cherkaev \cite{Gibiansky:1987:MCE}, and Allaire \cite{Allaire:1994:ELP},
but also provides a sharp bound on all possible $(\BGs^0,\BGe^0)$ pairs. Specifically,
given any pair $(\BGs^0,\BGe^0)$ that satisfies \eq{3.1} as an equality, there exists
a sequence of materials parameterized by $\Gd\to 0$, with porosity $1-f$,
and with effective elasticity tensors $\BC^*_\Gd$, such that
\beq \ajj{\BGs^0_\Gd=\BC^*_\Gd\BGe^0\to\BGs^0.}
\eeq{3.2}
\end{theorem}

\begin{theorem}
For $3d$-materials constructed from phase 1 with elasticity tensor $\BC_1$ defined by \eq{1.7a}
and a rigid material, with phase 1 occupying a volume fraction $f$, the complementary energy bound
\beq \widetilde{W}_f(\BGe^0)\leq\BGs^0:\BGe^0, \eeq{3.3}
where $\widetilde{W}_f(\BGe^0)$ can be computed from \eq{2.7} and \eq{2.8},
not only provides a sharp lower bound on the elastic energy, as established 
by Allaire and Kohn \cite{Allaire:1993:OBE},
but also provides a sharp bound on all possible $(\BGs^0,\BGe^0)$ pairs. Specifically,
given any pair $(\BGs^0,\BGe^0)$ that satisfies \eq{3.3} as an equality, there exists
a sequence of materials parameterized by $\Gd\to \infty$, with $f$ as the volume fraction of phase 1,
and with effective elasticity tensors $\BC^*_\Gd$, such that
\beq \BGe^0_\Gd=(\BC^*_\Gd)^{-1}\BGs^0\to\BGe^0.
\eeq{3.4}
\end{theorem}

In these theorems we assumed that $(\BGs^0,\BGe^0)$ were such that \eq{3.1} or \eq{3.3} held as an equality.
The case where one wants to achieve a pair $(\BGs^0,\BGe^0)$ such that \eq{3.1} or \eq{3.3} holds as an
inequality is then easily treated. Thus, for example, suppose one has inequality in \eq{3.1}. By taking extremely
thin struts in the pentamode construction in Figure 9 in \cite{Milton:2016:PEE} it is clear that for a given
$\BGs^0$ one can find a pentamode with elasticity tensor $\BC^*$ that supports the stress  $\BGs^0$
but which has a very large value of the compliance energy $\BGs^0:(\BC^*)^{-1}\BGs^0$. Laminating this
pentamode with a near optimal pentamode that supports the stress  $\BGs^0$ and almost achieves equality in \eq{3.1}, then
gives a pentamode material supporting the stress $\BGs^0$ that achieves the desired value of $\BGs^0:(\BC^*)^{-1}\BGs^0$.
Similarly, if one wants to achieve a pair $(\BGs^0,\BGe^0)$ such that \eq{3.3} holds as an inequality one takes
a unimode materials that has $\BGe^0$ as its easy mode of deformation, but has a very large value
of the elastic energy $\BGe^0:\BC^*\BGe^0$. Then one laminates this with the near optimal unimode that has the same
strain $\BGe^0$ as its easy mode of deformation to achieve a unimode material with $\BGe^0$ as its easy mode of deformation and
having the desired value of $\BGe^0:\BC^*\BGe^0$.

As mentioned in the introduction, this completely solves the weak $G$-closure problem for
$3d$-printed materials formed from one elastic material plus void,
and for materials that are composites of one elastic phase and a rigid phase. We assumed the elastic phase was 
isotropic with positive Lame modulus $\Gl$, but we could have allowed $\Gl$
to be negative or $\BC_1$
to be anisotropic, at some increased computational cost of evaluating
the functions $W_f(\BGs^0)$ and $\widetilde{W}_f(\BGe^0)$ according to the prescription
provided by Avellaneda \cite{Avellaneda:1987:OBM}. With the  weak $G$-closure problem
solved one can now formulate many optimization problems as relaxed problems and
proceed to solve them. However, if one wants to use the solution
of a relaxed problem to build realistic structures that approximately solve the 
original optimization problem one will need to search for geometries that have approximately
the same performance as the  multiscale near optimal pentamodes or multiscale near optimal unimodes,
whose construction was detailed
in \cite{Milton:2016:PEE} and \cite{Milton:2016:TCC}, respectively.

We remark that for $3d$-printed materials, near optimal pentamodes are not the only structures
that solve the weak $G$-closure problem. Indeed, if one has an elasticity tensor $\BC^*$ realizing a desired pair
$(\BGs^0,\BGe^0)$ then we are free to use a stiffer modified geometry, that is not necessarily a pentamode,
 with effective tensor $\underline{\BC}^*\geq \BC^*$
so long as $(\underline{\BC}^*-\BC^*)\BGe^0=0$. Similarly, for mixtures of the elastic phase with a rigid
phase, near optimal unimodes are not the only structures that solve the weak $G$-closure problem. 

%%%%%%%%%%%%%%%%%%%%%%%%%%%%%%%%%%%%%%%%%%%%%%%%%%%%%%%%%%%%%%%%%%%%%%%%%
\section*{Acknowledgements}
The authors thank the National Science Foundation for support through grant DMS-1211359. \ajj{Helpful feedback
from the referee is also gratefully acknowledged.}
%%%%%%%%%%%%%%%%%%%%%%%%%%%%%%%%%%%%%%%%%%%%%%%%%%%%%%%%%%%%%%%%%%%%%%%%%%%%%%%%%%%%%%%%%%

%\bibliographystyle{/home/milton/latex/camar/elsarticle/elsarticle-num.bst}
%\bibliographystyle{plain}
%\bibliography{/home/milton/tcbook,/home/milton/newref}
\end{document}